\documentclass[12pt,14paper]{article}
\usepackage{t1enc}
\usepackage[cp1251]{inputenc}
\usepackage[english]{babel}
\usepackage{graphicx}
\usepackage{epstopdf}
\usepackage{amsmath}
\usepackage{amssymb}
\usepackage{comment}

\newtheorem{theorem}{Theorem}

\newtheorem{lemma}{Lemma}

\def\square{\hbox{\vrule\vbox{\hrule\phantom{o}\hrule}\vrule}}

\def\Re{\mathop{\rm Re}\nolimits}
\def\Im{\mathop{\rm Im}\nolimits}

\def\diag{\mathop{\rm diag}\nolimits}

\newcommand{\trn}{^{\scriptscriptstyle {\rm T}}}

\begin{document}
\thispagestyle{empty}

\author{Ibragim A. Junussov
{\thanks{I.A. Junussov is with Department of Control Engineering,
Czech Technical University in Prague, e-mail: \nobreak{dxdtfxut@gmail.com}
The work was supported by the Ministry
of Education, Youth and Sports of the Czech \nobreak{Republic} within the project no.
CZ.1.07/2.3.00/30.0034 (Support for improving R\&D teams
and the development of intersectoral mobility at CTU in Prague). Author would
like to thank {A.L. Fradkov} for problem statement, Z. Hurak, A. Selivanov and \nobreak{I.~Herman} for useful discussions.
}}}
\title{Static consensus in passifiable linear networks}
\maketitle
\abstract{Sufficient conditions of consensus (synchronization) in networks described by digraphs and consisting
of identical determenistic SIMO systems are derived. 
Identical and nonidentical control gains (positive arc weights) are considered. Connection between 
admissible digraphs and nonsmooth hypersurfaces (sufficient gain boundary) is established. Necessary and sufficient conditions 
for static consensus by output feedback in networks consisting of certain class of double integrators
are rediscovered. Scalability for circle digraph in terms of 
gain magnitudes is studied. Examples and results of numerical simulations are presented.}
\section{Introduction}
Control of multi-agent systems has attracted significant interest in last decade since it has
a great technical importance~\cite{Bullo,mesb_eger_2010,saber_fax_murray07,ren_beard_acc05} and relates to
biological systems~\cite{reynolds_1987}.

In consensus problems agents communicate via decentralized controllers using relative measurements
with a final goal to achieve common behaviour (synchronization) which can evolve in time. 
Many approaches 
have been developed for a different problem settings.

Laplace matrix, its spectrum, and eigenspace plays crucial role in description and analysis of consensus 
problems. It has broad applications, e.g.~\cite{gutman_vidovic2002}. Not all possible digraph topologies
can provide consensus over dynamical networks. Admissible digraph topologies and connection with algebraic properties
of Laplace matrix have been found in~\cite{agaev_chebotarev_ait2000_en}. Analysis of tree strucure and Laplace matrices 
spectrum of digraphs are also studied by these authors. Work~\cite{agaev_chebotarev2013} contains examples of 
out-forests as well as useful graph theoretical concepts
and can be recommended as an entry reading to the research of these authors on algebraic digraph 
theory and consensus problems.

Concept of synchronization region in complex plane for a networks
consisting of linear dynamical systems is introduced in~\cite{Chen10}. In~\cite{Zhang_Lewis_Tac2011} 
this concept is used for analysis of synchronization with leader. Problem is solved using 
Linear Quadratic Regulator approach in cases when full state is available for measurement and 
when its not. In last case observers are constructed. 

Analysis of consensus with scalar coupling strenghts~\cite{Chen10,Zhang_Lewis_Tac2011}
is fruitful in a sense that conditions on gains (which depend on connection topology and
single agent properties) give more insight to problem. A lot of works on topic consider
dynamic couplings, however, for certain type of connections it might happen that
tunable parameters will exceed upper bound on possible control gains, i.e. won't meet physical limitations.
So, necessary and sufficient conditions on consensus achievement for different connection types in terms of 
coupling strenghts are needed. 

Celebrated Kalman-Yakubovich-Popov Lemma (Positive Real Lemma) establishes important connection between passivity
(positive-realness) of transfer function $\chi(s)$ and matrix relations on its minimal state-space realization
$(A,B,C),$ see~\cite{brogliato_dissipative_systems_book,kokotovic_arcak_autom2001}. Positive Real Lemma is a basis for Passification 
Method~\cite{Fradkov76_eng,F03EJC} (``Feedback Kalman-Yakubovich Lemma'') which answers question when a 
linear system can be be made passive, i.e. strictly 
positive-real (SPR) by static output feedback. Powerful idea of rendering system into passive by feedback
have been also studied for nonlinear systems, 
e.g.~\cite{byrnes_isidori_willems_1991,kokotovic_arcak_autom2001,kokotovic_1989}.

In consensus-type problems considering SPR agents with stable (Hurwitz) matrix $A$ leads to a synchronous
behaviour when all states going to zero. Latter is undesirable in essence, since such behaviour can be reached by 
local control without communication. So, instead of SPR systems it is possible to consider passifiable systems,
with an opportunity that a study is extendable to nonlinear systems. Also, Passification Method allows
to avoid constructing observers for reaching full-state consensus by output feedback. 
Observers implementation increase dimension of overall phase space and Complexity of a dynamic network.

In this paper Passification Method are used to synthesize a 
decentralized control law and to derive sufficient conditions of full state synchronization by relative 
output feedbacks in a networks described by digraphs with Linear Time Invariant dynamical nodes in continuous time. 
Assumptions made on network topology are minimal. Synchronous behaviour is described, including
case of nonidentical gains. It is determined that boundary of sufficient gain region geometrically is a 
hypersurface in corresponding gain space. For certain three node network this geometrical observation
connects algebraic properties of Laplace matrix with constructed hypersurface. Namely, Jordan block appears in
a direction of a cusp (nonsmooth) extremal point of the hypersurface. 

Necessary and sufficient conditions
for static consensus by output feedback in networks consisting of certain class of double integrators
have been rediscovered. Conditions are given in terms of Laplace matrix spectrum.

Scalability in a circle digraphs in terms of gain (coupling strength) is studied.
It is shown that common gain in large
cycle digraphs consisting of double integrators should grow not slower than quadratically in number of agents. 

Results of numerical simulations in 3 and 20 node double-integrator networks are presented.

\section{Theoretical study}
\subsection{Preliminaries and notations}
\label{pre}

Notations, some terms of graph theory and Passification Lemma are listed in this section.

\subsubsection{Notations}
Notation $\|.\|_2$ stands for Euclidian norm. For two symmetric matrices $M_1,M_2$ inequality
$M_1>M_2$ means that matrix $M_1-M_2$ is positive definite. Notation $\mathrm{col}(v_1,...,v_d)$
stands for vector $(v_1,...,v_d)\trn.$ Identity matrix of size $d$ is denoted by $I_d.$   
Vector $\mathbf{1}_d=(1,1,\ldots,1)$ is vector of size $d$ and consisting of ones.
Vector $\mathbf{0}_d$ is defined similarly. Matrix $\mathrm{diag}(v_1,\ldots,v_d)$ is square matrix
whose $i$-th element on main diagonal is $v_i, i=1,\ldots,d;$ other entries are zeroes.
Notation $\otimes$ stands for Kronecker product of matrices. Definition and properties of Kronecker
product, including eigenvalues property, can be found in ~\cite{bellman_eng,marcus_eng}. Direct sum
of matrices~\cite{horn_johnson_matr_analys} is denoted by $\oplus.$

\subsubsection{Terms of graph theory}
A pair $\mathcal{G}=(\mathcal{V},\mathcal{E}),$ 
where $\mathcal{V}$ -- set of vertices,
$\mathcal{E}\subseteq\mathcal{V}\times\mathcal{V}$ -- set of arcs (ordered pairs), is called digraph (directed graph). 
Let $\mathcal{V}$ have $N$  elements, $N\in\mathbb{N}.$ It is assumed hereafter
that graphs does not have self-loops, i.e. for any vertex $\alpha\in\mathcal{V}$ arc
$(\alpha,\alpha)\notin\mathcal{E}.$  

Digraph is called directed tree if all it vertices except one (called root) have exactly one parent
Let us agree that in any arc $(\alpha,\beta)\in\mathcal{E}$ vertex $\beta$ is parent or neighbour.
Directed spanning tree of a digraph $\mathcal{G}$ is a directed tree formed of all digraph $\mathcal{G}$ vertices and 
some of its arcs such that there exists path from any vertex to the root vertex in this tree. Existence of directed spanning tree 
has connection to principal achievement of synchronization in consensus-like problems.

A digraph is called weighted if to any pair of vertices $\alpha,\beta\in\mathcal{E}$ number
$w(\alpha,\beta)\geq 0$ is assigned such that:
$$w(\alpha,\beta)> 0 \quad\text{if}\quad (\alpha,\beta)\in\mathcal{E} \quad\text{and}\quad w(\alpha,\beta)=0 \quad\text{if}\quad (\alpha,\beta)\notin\mathcal{E}.$$
A digraph in which all nonzero weights are equal to $1$ will be referred as unit weighted.

An adjacency matrix $\mathcal{A}(\mathcal{G})$ is 
$N\times N$ matrix whose $i-\mbox{th},j-\mbox{th}$ entry is equal to $w(\alpha_i,\alpha_j), i,j=1,\ldots,N.$ 

Laplace matrix of digraph
$\mathcal{G}$ is defined as follows:
$$L(\mathcal{G})=\mathrm{diag}\left(\mathcal{A}(\mathcal{G})\cdot\mathbf{1}_N\right)-\mathcal{A}(\mathcal{G}).$$
Matrix $L(\mathcal{G})$ always has zero eigenvalue with corresponding right eigenvector
$\mathbf{1}_N: L(\mathcal{G})\cdot\mathbf{1}_N=0\cdot\mathbf{1}_N.$
By construction and Gershgorin Circle Theorem
all nonzero eigenvalues of $L$ have positive real parts. Let us denote by 
$v(L)\in\mathbb{R}^N$ left eigenvector of $L$ which is corresponding to zero eigenvalue 
and scaled such that $v(L)\trn \cdot\mathbf{1}_N=1.$ It is known that vector $v(L)$ describes synchronous 
behaviour if reached. 

Suppose that a digraph has directed spanning tree. A set of digraph vertices is called Leading Set
(``basic bicomponent'' in terms of~\cite{agaev_chebotarev2013}) if
subdigraph constructed of them is strongly connected and no vertex in this set has neighbours in the rest part of digraph.
Nonzero components of $v(L)$ and only them correspond to vertices of Leading Set. Definition of basic bicomponent 
is wider and applicable for digraphs with no directed spanning trees.

For illustration, by~\cite{harary_graph_theory}, there are 16 different types of digraphs which can be constructed on 3 nodes.
12 of them have directed spanning tree, among these 5 digraphs have Leading Set with 3 nodes, 2 digraphs have Leading 
Set with 2 nodes, and 5 digraphs have Leading Set with 1 node.                                                                              

\subsubsection{Passification Lemma}

Problem of linear system passification is a problem of finding static linear output feedback which is
making initial system passive. It was solved in~\cite{Fradkov76_eng,F03EJC} for nonsquare SIMO and 
MIMO systems including case of complex parameters. Brief outline of SIMO systems passification is 
given below.  

Let $A,B,C$ be real matrices of sizes $n\times n, n\times 1, n\times l$ accordingly. Denote by
$\chi(s)=C\trn(sI-A)^{-1}B, s\in\mathbb{C}.$ Let vector $g\in\mathbb{R}^l.$ If numerator of
function $g\trn\chi(s)$ is Hurwitz with degree $n-1$ and has positive coefficients then  function 
$g\trn\chi(s)$ is called hyper-minimum-phase. 

\begin{lemma}
({\it Passification Lemma}~\cite{Fradkov76_eng,F03EJC})
Following statements are equivalent.
\begin{enumerate}
\item There exists vector $g\in\mathbb{R}^l$ such that function 
$g\trn\chi(s)$ is hyper-minimum-phase. 
\item Number $\varkappa_0=\sup\limits_{\omega\in\mathbb{R}^1}{\rm Re} \big
(g\trn \chi(\mathfrak{i}\omega)\big)^{-1}$ is positive $\varkappa_0>0$ and for any $\varkappa>\varkappa_0$
there exists $n\times n$ real matrix $H=H\trn >0$ satisfying following matrix relations
\begin{equation}
\label{h_m_p}
HA_*+A_*\trn H<0,\quad HB=Cg, \quad A_*=A-\varkappa Bg\trn C\trn.
\end{equation}
\end{enumerate}
\end{lemma}

\subsection{Problem statement and assumptions}
\label{probl_statement}
Consider a network consisting of $N$ agents modelled as linear dynamical systems:
\begin{equation}
\label{main_system}
\begin{aligned}
\dot{x}_i(t) &= Ax_i(t)+Bu_i(t),\\
y_i(t) &= C\trn x_i(t), 
\end{aligned}
\end{equation}
where $i=1,\ldots,N, x_i\in\mathbb{R}^n$ -- state vector, $y_i\in\mathbb{R}^l$ -- output 
or measurements vector, $u_i\in\mathbb{R}^1$ -- input or control, $A,B,C$ are real matrices 
of according size. By associating agents with $N$ vertices of unit weighted digraph $\mathcal{G}$ and
introducing set of arcs one can
describe information flow in the network.
For $i=1,\ldots,N$ let us introduce notation for relative outputs 
$$\overline{y}_i(t)=\sum_{j\in\mathcal{N}_i}(y_i(t)-y_j(t)),$$
where $\mathcal{N}_i$ is a set of $i$-th agents neighbours.

Problem is to design controllers which use relative outputs and ensure achievement of
the state synchronization (consensus) of all agents:
\begin{equation}
\label{goal}
\lim_{t\to\infty}(x_i(t)-x_j(t))=0, \quad i,j=1,\ldots,N.
\end{equation}
In the case of synchronization achievement asymptotical behaviour of all agents will be 
described by same time-dependant consensus vector which is denoted hereafter by $c(t):$
$$
\lim_{t\to\infty}(x_i(t)-c(t))=0, \quad i=1,\ldots,N.
$$
Let us make following assumption about dynamics of a single agent.

{A1)} \textit{There exists vector $g\in\mathbb{R}^l$ such that transfer function 
$g\trn\chi(s)=g\trn C\trn (sI_n-A)^{-1}B$ is hyper-minimum-phase.}



Now let us make an assumption on graph topology.

{A2)} \textit{Digraph $\mathcal{G}$ has at least one directed spanning tree.}

Zero eigenvalue of Laplace matrix $L$ has unit multiplicity iff this assumption holds~\cite{agaev_chebotarev_ait2000_en}.

\subsection{Static identical control}
\label{identical_contr}
Denote $r(L)=\min\limits_{\lambda_i\neq 0} \Re \lambda_i$ where $\lambda_i$ are 
eigenvalues of $L.$ Under assumption A2 zero eigenvalue is simple. By properties of $L$ 
other eigenvalues lie in open right half of complex plane, so $r(L)$ is positive number. 

Suppose that assumption A1 holds with known vector $g\in\mathbb{R}^l.$
Consider following static consensus controller with gain  
$k\in\mathbb{R}^1, k>0$ which is same for all agents:
\begin{equation}
\label{stat_contr}
u_i(t)=-k g\trn \overline{y}_i(t),\quad i=1,\ldots,N,
\end{equation}
where relative output $\overline{y}_i(t)$ has been defined in previous section. 
Denote $x(t)=\mathrm{col}(x_1(t),\ldots,x_N(t)).$
\begin{theorem}
\label{th1}
Let assumptions A1 and A2 hold. Then for all $k$ such that
\begin{equation}
\label{k}
k>\frac{\varkappa_0}{r(L)}
\end{equation}
controller \eqref{stat_contr} ensures achievement of goal \eqref{goal} in dynamical network \eqref{main_system};
asymptotical behaviour in the case of goal achievement is described by following consensus vector
\begin{equation}
\label{c_t_ident}
c(t)=\exp(At)(v(L)\trn\otimes I_n)x(0).
\end{equation}
\end{theorem}

\textit{Proof.} Closed loop system \eqref{main_system}, \eqref{stat_contr} can be rewritten in a following form
\begin{equation}
\label{closed_lp_st}
\dot{x}(t)=(I_N\otimes A - k L\otimes B g\trn C\trn)x(t).
\end{equation}
Consider nonsingular matrix $P$ (real or complex) such that
\begin{equation*}
\Lambda=
\begin{pmatrix}
0     &\mathbf{0}_{N-1}\trn\\
\mathbf{0}_{N-1} &\Lambda_e\\
\end{pmatrix}=P^{-1}LP,
\end{equation*}
where $\Lambda_e\in\mathbb{R}^{(N-1)\times(N-1)}$ or $\Lambda_e\in\mathbb{C}^{(N-1)\times(N-1)}.$ All eigenvalues of $\Lambda_e$ have positive 
real parts. By considering first (zero) columns of matrices $P\Lambda=LP$ and 
$(P\trn)^{-1}\Lambda\trn=L\trn (P^{-1})\trn$ we can accept that first column of $P$ is $\mathbf{1}_N$
and first row of $P^{-1}$ is $v(L)\trn.$

Let us apply coordinate transformation $z(t)=(P^{-1}\otimes I_n)x(t)$ and rewrite \eqref{closed_lp_st}:
\begin{equation}
\label{z_1}
\dot{z}_1(t)=Az_1(t),
\end{equation}
\begin{equation}
\label{z_e}
\dot{z}_e(t)=\left((I_{N-1}\otimes A)-k(\Lambda_e\otimes Bg\trn C\trn)\right)z_e(t),
\end{equation}
where $z=\mathrm{col}(z_1,z_e), z_1\in\mathbb{R}^n$ or $z_1\in\mathbb{C}^n.$ 
Note that zero solution of \eqref{z_e} is globally asymptotically stable iff
goal \eqref{goal} is achieved.

For simplicity let $P, \Lambda_e,$ and $z(t)$ be real till the end of proof.
For any fixed $k$ satisfying \eqref{k} there exists $0<\varepsilon_s<1$ such that 
$$
\varepsilon_s k>\frac{\varkappa_0}{r(L)}.
$$

Eigenvalues of matrix $(\Lambda_e-\varepsilon_s r(L)I_{N-1})$ have positive real parts.
Therefore, according to~\cite{lmi_book}, there exists $(N-1)\times(N-1)$ real 
matrix $Q=Q\trn>0$ such that following Lyapunov inequality holds
$$
(\Lambda_e-\varepsilon_s r(L)I_{N-1})\trn Q + Q (\Lambda_e-\varepsilon_s r(L)I_{N-1}) >0.
$$

We can rewrite last inequality
$$
\Lambda_e\trn Q + Q \Lambda_e > 2 \varepsilon_s r(L) Q.
$$

By assumption A1 there exists $H=H\trn>0$ such that \eqref{h_m_p} is true with 
$\varkappa= \varepsilon_s k r(L),$ since $\varkappa>\varkappa_0.$
Considering following Lyapunov function 
$$
V(z_e(t))=z_e\trn(t)(Q\otimes H)z_e(t)
$$
and derivativing it along the nonzero trajectories of \eqref{z_e}, we obtain
\begin{equation*}
\begin{aligned}
\frac{\mathrm{d}}{\mathrm{dt}}V(z_e(t)) &=z_e\trn(t)(Q\otimes (A\trn H+HA)-k(\Lambda_e\trn Q+Q\Lambda_e)\otimes (Cgg\trn C\trn))z_e(t)\leq\\
&\leq z_e\trn(t)(Q\otimes(A\trn H+ HA)-2k\varepsilon_s r(L) Q\otimes (Cgg\trn C\trn))z_e(t)=\\
&= z_e\trn(t)(Q\otimes((A\trn-\varkappa CgB\trn)H+H(A-\varkappa Bg\trn C\trn)))z_e(t)=\\
&= z_e\trn(t)(Q\otimes(A_*\trn H+HA_*))z_e(t)<0.
\end{aligned}
\end{equation*}
Matrix relations~\eqref{h_m_p} have been used here. Last inequality concludes the proof.\square

Assumptions of this Theorem are relaxed in comparison with Theorem 2 from~\cite{fradkov_junussov_cdc11}.
Proof of Theorem~\ref{th1} also provides following auxiliary result.

\begin{lemma}
\label{lemma0}
Let assumption A2 hold. 
Controller \eqref{stat_contr} ensures achievement of goal \eqref{goal} in dynamical network \eqref{main_system}
if, and only if, all eigenvalues of matrix
$$
R=(I_{N-1}\otimes A)-k(\Lambda_e\otimes Bg\trn C\trn)
$$
have negative real parts. In the case of goal achievement 
asymptotical behaviour is described by \eqref{c_t_ident}.
\end{lemma}

\subsection{Nonidentical control and Gain Region}
\label{gain_region}
Let the initial digraph $\mathcal{G}$ be unit weighted.
Let us fix Laplace matrix $L$ and consider static control with nonidentical gains $k_i>0:$
\begin{equation}
\label{gain_nonident}
u_i(t)=-k_i g\trn \overline{y}_i(t),\quad i=1,\ldots,N.
\end{equation}
Without loss of generality we can assume that
network does not have a leader (formally: cardinality of Leading Set is more than 1), since in leader case we can reduce following 
consideration of synchronization gain region to lower dimension $N-1.$ 

Let us denote by $\hat{k}=(k_1,\ldots,k_N)$ and by $\hat{k}'=(k_1',\ldots,k_N')$ point which is projection of
point $\hat{k}$ on unit sphere $\mathcal{S}:$
$$
\hat{k}=k\cdot\hat{k}',\quad k>0,\quad \sum_{i=1}^N(k_i')^2=1,
$$
where scalar common gain $k$ is radius vector magnitude of point $\hat{k}.$
Points $\hat{k},\hat{k}'$ lie in orthant 
$\mathcal{O}=\{(k_1,\ldots,k_N)\in\mathbb{R}^N| k_i>0,\,\,i=1,\ldots,N\}.$

Denote by $K'=\mathrm{diag}(k_1',k_2',\ldots,k_N').$ Laplace matrices $L$ and $K'L$ correspond to
same digraphs which differ only in arc weights. Equation
for closed loop system \eqref{main_system},\eqref{gain_nonident} can be rewritten as follows
$$
\dot{x}(t)=(I_N\otimes A - k (K'L\otimes B g\trn C\trn))x(t).
$$
By repeating proof of Theorem \ref{th1} we can formulate following result.
\begin{theorem}
\label{th2}
Let assumptions A1 and A2 hold. Then for all $k_i=k\cdot k_i'$ such that
$$
\sum_{i=1}^N(k_i')^2=1,\quad k>\frac{\varkappa_0}{r(K'L)}
$$
controller \eqref{gain_nonident} ensures achievement of goal \eqref{goal} in dynamical network \eqref{main_system};
asymptotical behaviour in the case of goal achievement is described by following consensus vector
\begin{equation}
\label{cons_vect}
c(t)=\exp(At)(v(K'L)\trn\otimes I_n)x(0).
\end{equation}
\end{theorem}

Denote by $\mathcal{K}\subset\mathcal{O}$ region in orthant such that for any $(k_1,\ldots,k_N)\in\mathcal{K}$
control \eqref{gain_nonident} ensures achievement of the goal \eqref{goal} in network \eqref{main_system}, 
\eqref{gain_nonident}. Consider following region
$$
\mathcal{K}_r=\left\{\hat{k}\in\mathcal{O}\left|\hat{k}=k\cdot\hat{k}',\hat{k}'\in\mathcal{S},k>\frac{\varkappa_0}{r(K'L)}\right.\right\}
$$
which is subset of $\mathcal{K}: \mathcal{K}_r\subset\mathcal{K}.$
Let us consider closed part of unit sphere
$$\mathcal{S}_\varepsilon=\left\{\hat{k}'\in\mathcal{S}\left|k_i'\geq\varepsilon,i=1,\ldots,N\right.\right\}, 
\quad\varepsilon>0.$$
Point on $\mathcal{S}_\varepsilon$ determines ray (half-line) in $\mathcal{O}$ with initial point at the origin. 
According to Theorem \ref{th2}, by moving along this ray from origin, i.e. increasing $k,$ we 
will reach $\mathcal{K}_r.$ 
Consider map
$$
h:k'\mapsto\,\frac{\varkappa_0}{r(K'L)}\cdot k',\quad k'\in \mathcal{S}_\varepsilon,
$$
which is continuous as a composition of continuous maps (~\cite{horn_johnson_matr_analys}, continuous dependence 
of matrix eigenvalues on parameters). Image of this map is a subset of boundary $\partial\mathcal{K}_r,$ therefore,
by continuity of map $h,$ boundary $\partial\mathcal{K}_r$ is a hypersurface in $\mathbb{R}^N.$
Further, let us consider induced map
$$
h_\rho:k'\mapsto\, \| h(k') \|_2,\quad k'\in \mathcal{S}_\varepsilon.
$$

Domain
$\mathcal{S}_\varepsilon$ is compact, so we can apply Weierstrass Extreme Value Theorem and arrive at following lemma.
\begin{lemma}
\label{lemma1}
Map $h:\mathcal{S}_\varepsilon\to h(\mathcal{S}_\varepsilon)\subset\partial\mathcal{K}_r$ is continuous.
Map $h_\rho:\mathcal{S}_\varepsilon\to \mathbb{R}^1$ is continuous and has minimum and maximum.
\end{lemma}
Generally, hypersurface $\partial\mathcal{K}_r$ is not smooth in all its points. Alternatively, part of a simplex
of according dimension can be taken instead of sphere part to serve as the domain for maps $h$ and $h_\rho.$

Pairwise ratios of nonidentical gains and common gain define homogeneous coordinates in orthant. Common gain $k$ coefficient
relates to reachability of consensus and to speed of convergence but it doesn't influence consensus vector. Also,
consensus vector can be changed only 
by gain ratios variation within Leading Set of agents, see~\cite{agaev_chebotarev_cycl2012}.


\section{Double-integrator networks}

\subsection{Agents description}
\label{agents_descr}
Suppose that each agent $S_i$ in a network is modelled as follows
\begin{equation}
\label{double_integr}
\begin{aligned}
&\dot{x}_i=Ax_i+Bu_i, \qquad y_i=C\trn x_i,\quad i=1,\ldots,N,\\
&A=
\begin{pmatrix}
0&0\\
1&0\\
\end{pmatrix},\quad
B=
\begin{pmatrix}
2\\
0\\
\end{pmatrix},\quad
C=
\begin{pmatrix}
0.5\\
0.5\\
\end{pmatrix}.
\end{aligned}
\end{equation}
For $g=1$  transfer function $g\trn\chi(s)=C\trn(sI_2-A)^{-1}B=\frac{s+1}{s^2}$ is hyper-minimum-phase.
It can be shown that number $\varkappa_0=1.$ 

First and second components of ${x}_i$ can describe (or can be interpreted as) velocity and position.
Single system \eqref{double_integr} can be viewed as double integrator with transfer function $1/s^2$ and
proportionally differential (PD) control applied to it.

Since $g=1,$ static consensus controller \eqref{stat_contr} has following form:
\begin{equation}
\label{stat_contr_di}
u_i(t)=-k \overline{y}_i(t),\quad i=1,\ldots,N.
\end{equation}

\subsection{Necessary and sufficient conditions on consensus}
\label{cycle_digraph}
Let us denote by $L_N^C$ Laplace matrix of unit weighted cycle digraph which is
consisting of $N$ nodes $S_j$ with exactly $N$ arcs 
$$(S_1,S_2)\cup\ldots\cup(S_{j},S_{j+1})\cup\ldots\cup(S_{N-1},S_N)\cup(S_N,S_1).$$
Eigenvalues of $L_N^C$ are evenly located at circle in complex plane~\cite{FaxMurray04}:
\begin{equation}
\label{lambda_circle}
\lambda_j=1-\exp{\left(\mathfrak{i}\cdot j\cdot\frac{2\pi}{N}\right)}, j=0,\ldots,N-1, \mathfrak{i}^2=-1.
\end{equation}

\begin{theorem}
\label{th_circle}
Controller \eqref{stat_contr_di} ensures achievement of goal \eqref{goal} in dynamical network consisting
of $N$ double integrators \eqref{double_integr}
connected in directed cycle if, and only if,
\begin{equation}
\label{ctg_law}
k>\frac 12 \cot^2 \frac {\pi}N.
\end{equation}
\end{theorem}
\textit{Proof.} Let us diagonilize $L_N^C.$ Matrix $R$ from Lemma~\ref{lemma0} in our case is block diagonal
\begin{equation*}
R=R_1\oplus R_2\oplus\ldots\oplus R_{N-1},
\end{equation*}
where
\begin{equation*}
\begin{aligned}
&R_j=
\begin{pmatrix}
-k\lambda_j&-k\lambda_j\\
1&0\\
\end{pmatrix},\quad j=1,\ldots,N-1.
\end{aligned}
\end{equation*}
So, matrix $R$ is stable iff matrices $R_j$ are stable for all $j=1,\ldots,N-1.$
Characteristic polynomial of $R_j$ is 
\begin{equation}
\label{f_j}
f_j(z)=z^2+k\lambda_j z+k\lambda_j, j=1,\ldots,N-1.
\end{equation}
Let $k\lambda_j=\alpha_j+\mathfrak{i}\beta_j, \alpha_j,\beta_j\in\mathbb{R},j=1,\ldots,N-1.$
Taking in account~\eqref{lambda_circle} we can obtain 
$$\alpha_j=2k\sin^2 \frac{j\pi}{N},\quad \beta_j=-2k\left(\sin \frac{j\pi}{N}\right)\left(\cos \frac{j\pi}{N}\right),\quad j=1,\ldots,N-1.$$

Now let argument of $f_j(z)$ run on imaginary axis and let us decompose this polynomial on real and imaginary parts:
$$f_j(\mathfrak{i}\,\omega)=\varphi_j(\omega)+\mathfrak{i}\psi_j(\omega),\quad \omega\in\mathbb{R}^1,$$
where $j=1,\ldots,N-1$ and:
$$\varphi_j(\omega)=-\omega^2-\beta_j\cdot\omega+\alpha_j,\quad  \psi_j(\omega)=\alpha_j\cdot\omega+\beta_j.$$

According to Hermite-Biehler Theorem, polynomial $f_j(z)$ is stable iff both of following conditions satisfied:
\begin{itemize}
\item roots of $\varphi_j(\omega)$ and $\psi_j(\omega)$
are interlacing;
\item Wronskian is positive
$$\varphi_j(\omega_0)\cdot\psi^\prime_j(\omega_0)-\varphi^\prime_j(\omega_0)\cdot\psi_j(\omega_0)>0$$
for at least one value of argument $\omega_0.$
\end{itemize}
Wronskian is positive for $\omega_0=0,j=1,\ldots,N-1.$ Root interlacing property is equivalently transformable
to 
$$
k>\frac 12 \cot^2\frac{j\pi}{N}, j=1,...,N-1.
$$
Right parts of these $N-1$ inequalities reach maximum when $j=1$ (also when $j=N-1$) and this concludes proof.\square

Therefore, for a large increasing number of 
agents $N$ gain $k$ should grow as $N^2:$
\begin{equation}
\label{equiv}
k\sim \frac {N^2} {2\pi^2},\quad N\to\infty.
\end{equation}

It is possible to 
conclude that consensus in large cycle digraphs is hard to achieve, at least for agents \eqref{double_integr}, 
since an arbitrary high gains are not physically realizable.

On other hand, it is worth noting that cycle digraph is the graph with minimal number of edges which is delivering
average consensus among all its nodes, it is strongly connected. 

Remark (see~\cite{agaev_chebotarev_cycl2012}) \textit{Minimality in edges number provides with simple 
relations on  nonidentical gains and left eigenvector $v(KL^C_N)$ components for agents in form \eqref{main_system}:
\begin{equation}
\label{circle_relation}
k_1 v_1= k_2 v_2 =\ldots=k_N v_N,
\end{equation}
i.e. the less coupling strength agent have the more it impacts Synchronous Behaviour.
}

In other words, all pairs $(k_j,v_j)$ lie on same hyperbola.

Following result can be obtained by repeating proof of Theorem \eqref{th_circle}.
\begin{theorem}
\label{th_real_spectr}
Consider network $S$ consisting of $N$ agents \eqref{double_integr}. Let a digraph, describing information flow,
contain directed spanning tree. Let Laplace matrix of the digraph  have real spectrum.
Controller \eqref{stat_contr_di} ensures achievement of goal \eqref{goal} in dynamical network consisting
of $N$ double integrators \eqref{double_integr} if, and only if,
\begin{equation*}
k>0.
\end{equation*}
\end{theorem}
\textit{Proof.} For diagonalizable Laplace matrix with real spectrum statement is following from well-known fact that 
polynomial \eqref{f_j} with real coefficients is stable iff its coefficients are positive. For nondiagonalizable 
Laplace matrix $L$ let us transform it to Jordan form. Expansion of matrix $R-zI$ determinant
shows that only 
determinants $R_j-zI_2$ across main diagonal are forming (factorizing)
characteristic polynomial of $R.$ \square

Note that undirected graphs (i.e. digraphs with symmetric $L$)
have real spectrum and some class of digraphs have real spectrum too, 
e.g. directed path graphs~\cite{herman_platoon_2014}.

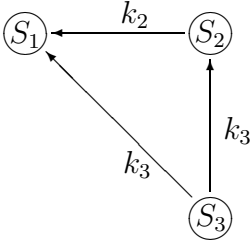
\begin{figure}
\centering
\begin{picture}(200,100)

\put(40,80){\circle{17}}
\put(34,76.4){$S_1$}

\put(110,80){\circle{17}}
\put(104,76.6){$S_2$}

\put(110,10){\circle{17}}
\put(104.4,7){$S_3$}
\put(100,80){\vector(-1,0){50}}
\put(110,20){\vector(0,1){50}}
\put(103,18){\vector(-1,1){55}}
\put(76,84){$k_2$}
\put(115,40){$k_3$}
\put(77,27){$k_3$}
\end{picture}
\caption{Digraph of 3 nodes}
\label{gr3nodes}
\end{figure}

We can formulate similar result for general digraphs.

\begin{theorem}
\label{th_compl_spectr}
Consider network $S$ consisting of $N$ agents \eqref{double_integr}. Let digraph $\mathcal{G}$, describing information flow,
contain at least one directed spanning tree. Let all nonzero eigenvalues of Laplace matrix $L(\mathcal{G})$ be denoted by 
$\lambda_j,1\leq j\leq N-1.$ Controller \eqref{stat_contr_di} ensures achievement of goal \eqref{goal} in dynamical network consisting
of $N$ double integrators \eqref{double_integr} if, and only if,
\begin{equation}
k>\max_{1\leq j\leq N-1}\frac {\sin^2 (\arg(\lambda_j))}{\Re{\lambda_j}}.
\end{equation}
\end{theorem}
\textit{Proof.} Using similar argumentation as in proofs of Theorems \ref{th_real_spectr} and  \ref{th_circle} 
we arrive at study of polynomials \eqref{f_j} stability with 
$$
\alpha_j=k\cdot\Re(\lambda_j), \beta_j=k\cdot\Im(\lambda_j), j=1,\ldots,N-1.
$$
Wronskian property does hold for $\omega_0=0.$ For $j=1,\ldots,N-1$ root interlacing property is equivalent to
trigonometric inequality
$$
\alpha_j+\beta_j\tan(\arg(\lambda_j))>\tan^2(\arg(\lambda_j))
$$
or
$$
k\Re(\lambda_j)(1+\tan^2(\arg(\lambda_j)))>\tan^2(\arg(\lambda_j)).
$$\square
\section{Examples and numerical simulations results}
\subsection{Three node digraph and gain region}
\label{3digraph}
\begin{figure}
\minipage{0.1\textwidth}
\endminipage\hfill
\minipage{0.4\textwidth}
\includegraphics[width=\linewidth]{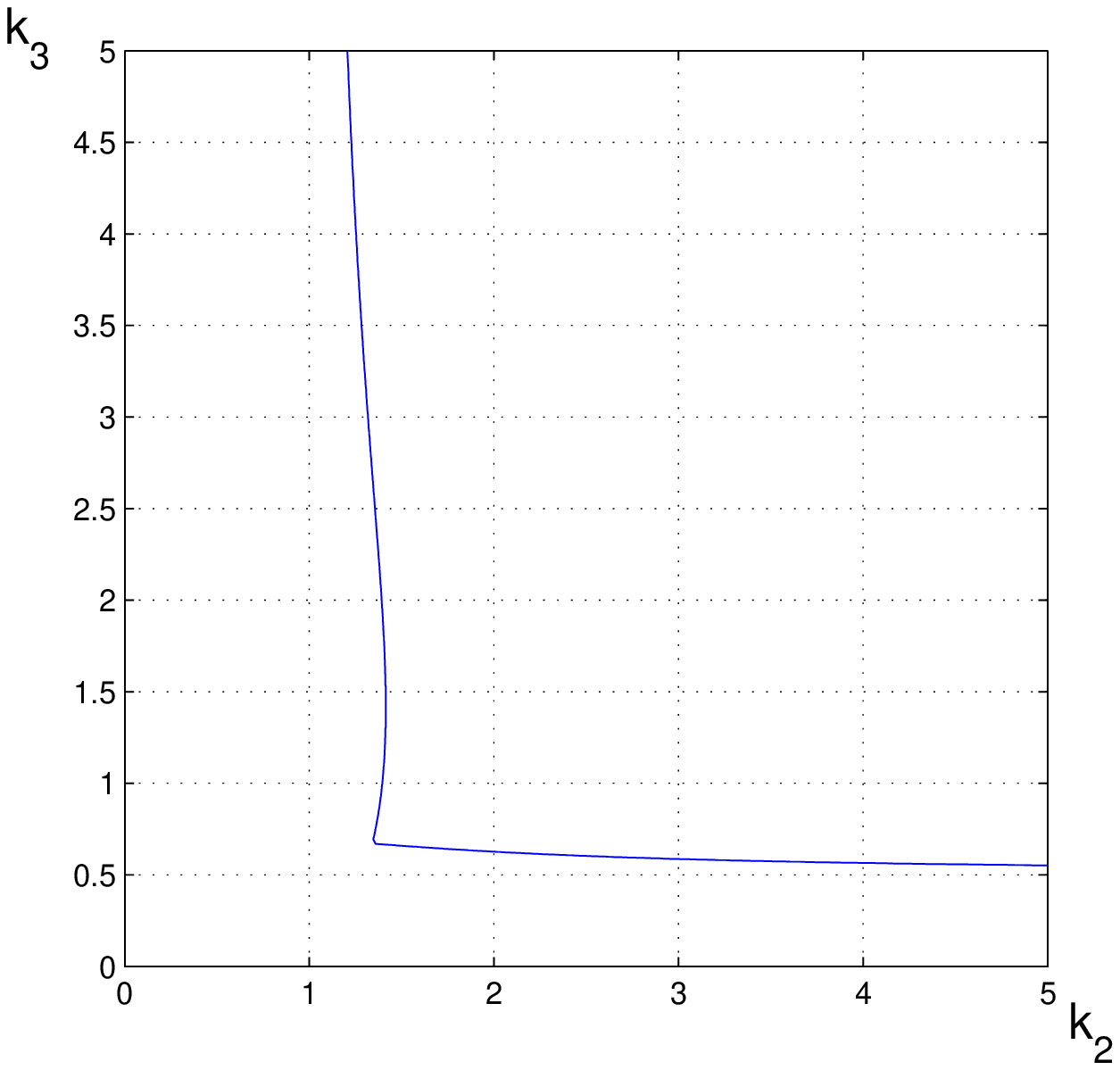}
\caption{Gain region boundary $\partial\mathcal{K}_r$.}
\label{boundary_dkr}
\endminipage\hfill
\minipage{0.4\textwidth}
\includegraphics[width=\linewidth]{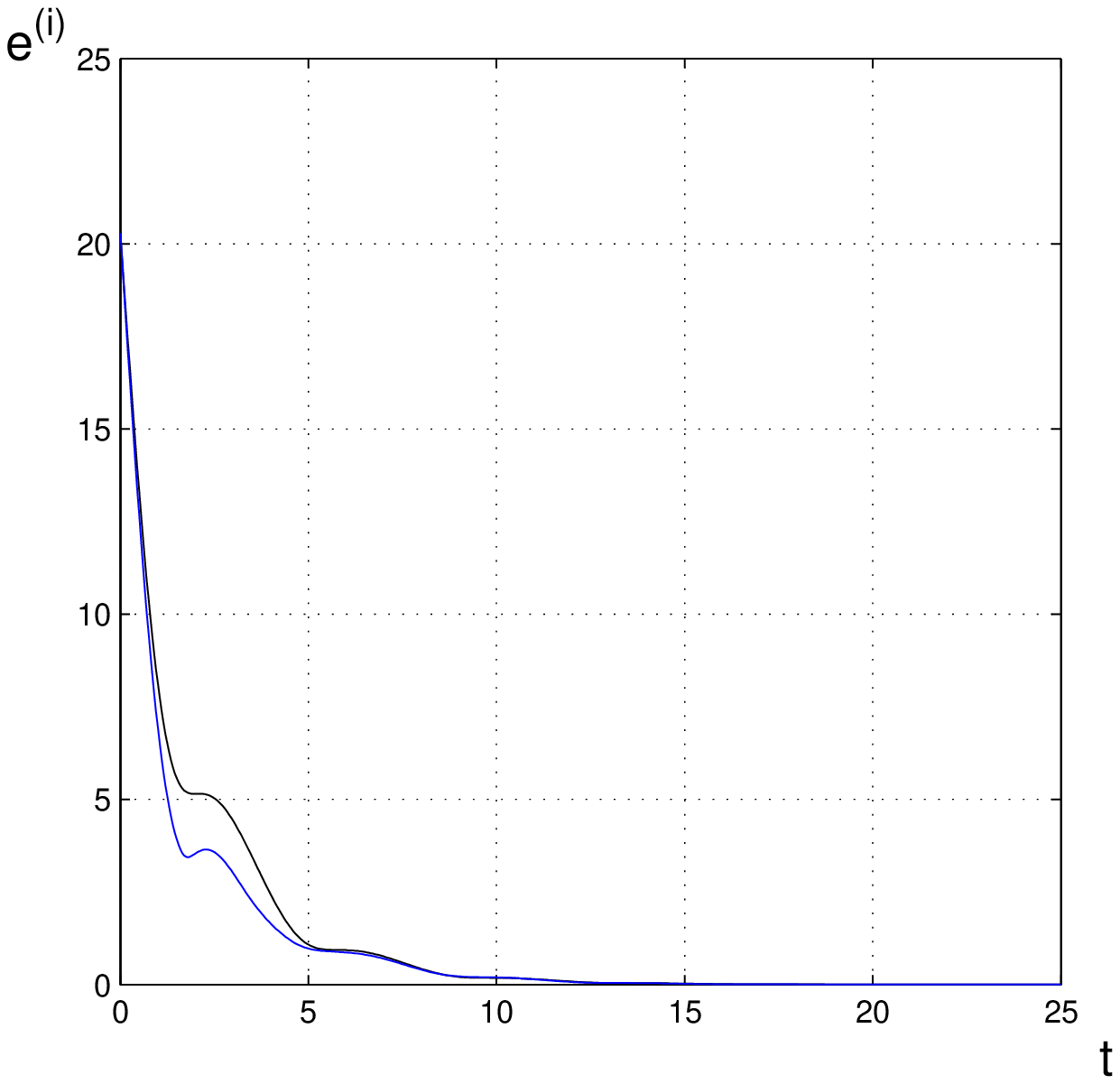}
\caption{Performance with identical ($e^{(1)}(t)$; black line) and nonidentical ($e^{(2)}(t)$; blue line) gains.}
\label{e1e2}
\endminipage\hfill
\minipage{0.4\textwidth}
\endminipage\hfill
\minipage{0.1\textwidth}
\endminipage\hfill
\end{figure}

Consider digraph shown on Fig.~\ref{gr3nodes} with dynamic nodes described in section~\ref{agents_descr}.
By Lemma ~\ref{lemma1} distance from origin to $\partial\mathcal{K}_r$ reaches minimum. Let us draw
$\partial\mathcal{K}_r.$ First, let $\varepsilon,\delta\in\mathbb{R}.$ Let $\varepsilon>0$ be small, 
and $\varepsilon\leq\delta\leq 1-\varepsilon, k_2'=\delta, k_3'=1-\delta.$ Let $L$ be unit weighted.
Eigenvalues of matrix $\mathrm{diag}(0,k_2',k_3')L$ are real: $\{0,\delta,2-2\delta\}.$ Using Theorem~\ref{th_real_spectr}
we conclude that $\mathcal{K}=\{k_2>0,k_3>0\}.$ Any $\delta\in[\varepsilon,1-\varepsilon]$ determines angle 
$$
\gamma(\delta)=\arctan\frac{1-\delta}\delta=\arctan\frac{k_3'}{k_2'}=\arctan\frac{k_3}{k_2},
$$
and radius vector $\rho(\delta)=1/\min_{\delta\in[\varepsilon,1-\varepsilon]}\{\delta,2-2\delta\}.$
Note that pair $(\rho(\delta),\gamma(\delta))$ is polar coordinates of boundary $\partial\mathcal{K}_r.$
We can conclude that minimum on $\rho(\delta)$ is realized on a point for which $\delta=2/3,$ and $k_2 : k_3 = 2 : 1.$
Boundary $\partial\mathcal{K}_r$ is presented on Fig.~\ref{boundary_dkr}. For all $(k_2,k_3)=k\cdot(2,1)$
matrix $K'L$ is similar to 
$$0\oplus \begin{pmatrix}
2k&1\\
0&2k\\
\end{pmatrix}.$$

Let us consider two cases:
\begin{enumerate}
\item $\delta=1/2,$ identical gains $k_2^{(1)}=k_3^{(1)}=0.527\cdot k;$
\item $\delta=2/3,$ nonidentical gains $k_2^{(2)}=\frac 23 \cdot k,\quad k_3^{(2)}=\frac 13 \cdot k.$
\end{enumerate}
By Theorem~\ref{th2} common gain $k$ is as follows
$k=3/2={\varkappa_0}/{r(K^{(2)}L)},\,
K^{(2)}=\diag(0,k^{(2)}_2,k^{(2)}_3).$
Identical gains are chosen such that $\left\|\left(k_2^{(1)},k_3^{(1)}\right)\right\|_2\approx\left\|\left(k_2^{(2)},k_3^{(2)}\right)\right\|_2.$ 

Let us choose following initial conditions:
$$
x_1(0)=\mathrm{col}(2,-2), 
x_2(0)=\mathrm{col}(-7,3), 
x_3(0)=\mathrm{col}(1,-3).
$$
Denote by $e(t)=\sum_{i=1}^2\|x_i(t)-x_{i+1}(t)\|_2$ sum of error norms or disagreement measure;
$e^{(1)}(t)$ error in the first case, $e^{(2)}(t)$ error in the second case. Results of 25 sec.
 simulations are shown on Fig.~\ref{e1e2}. 

Note that consensus vector \eqref{cons_vect} does not changes for all 
$(k_2, k_3)\in\mathcal{K}$ since subsystem $S_1$ is leader.

\subsection{Twenty node digraph and nonidentical control}
\label{dodec_digraph}
\begin{figure}
\minipage{0.1\textwidth}
\endminipage\hfill
\minipage{0.4\textwidth}
\includegraphics[width=8in]{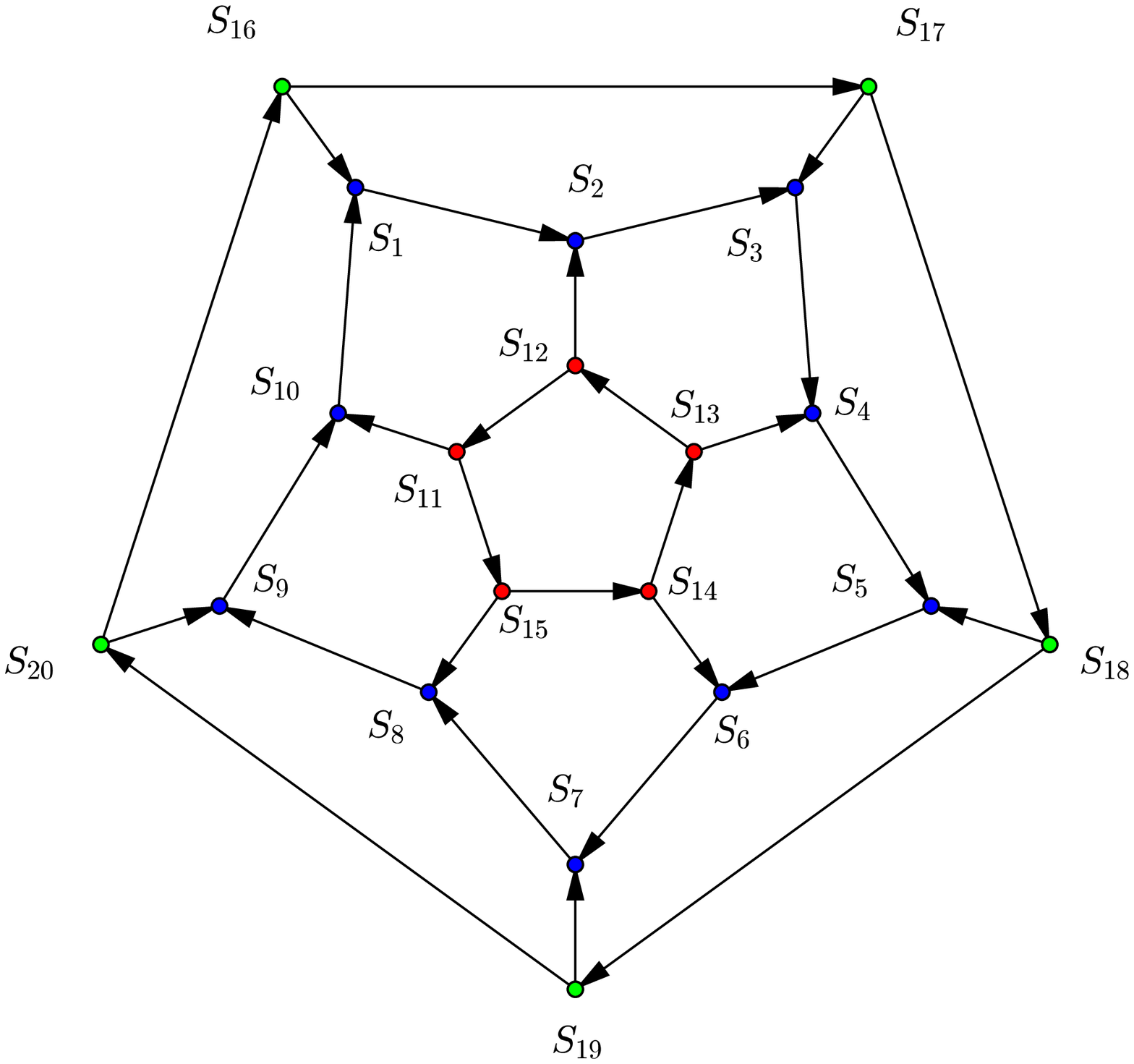}
\caption{Digraph of 20 agents.}
\label{dodec}
\endminipage\hfill
\minipage{0.4\textwidth}
\includegraphics[width=\linewidth]{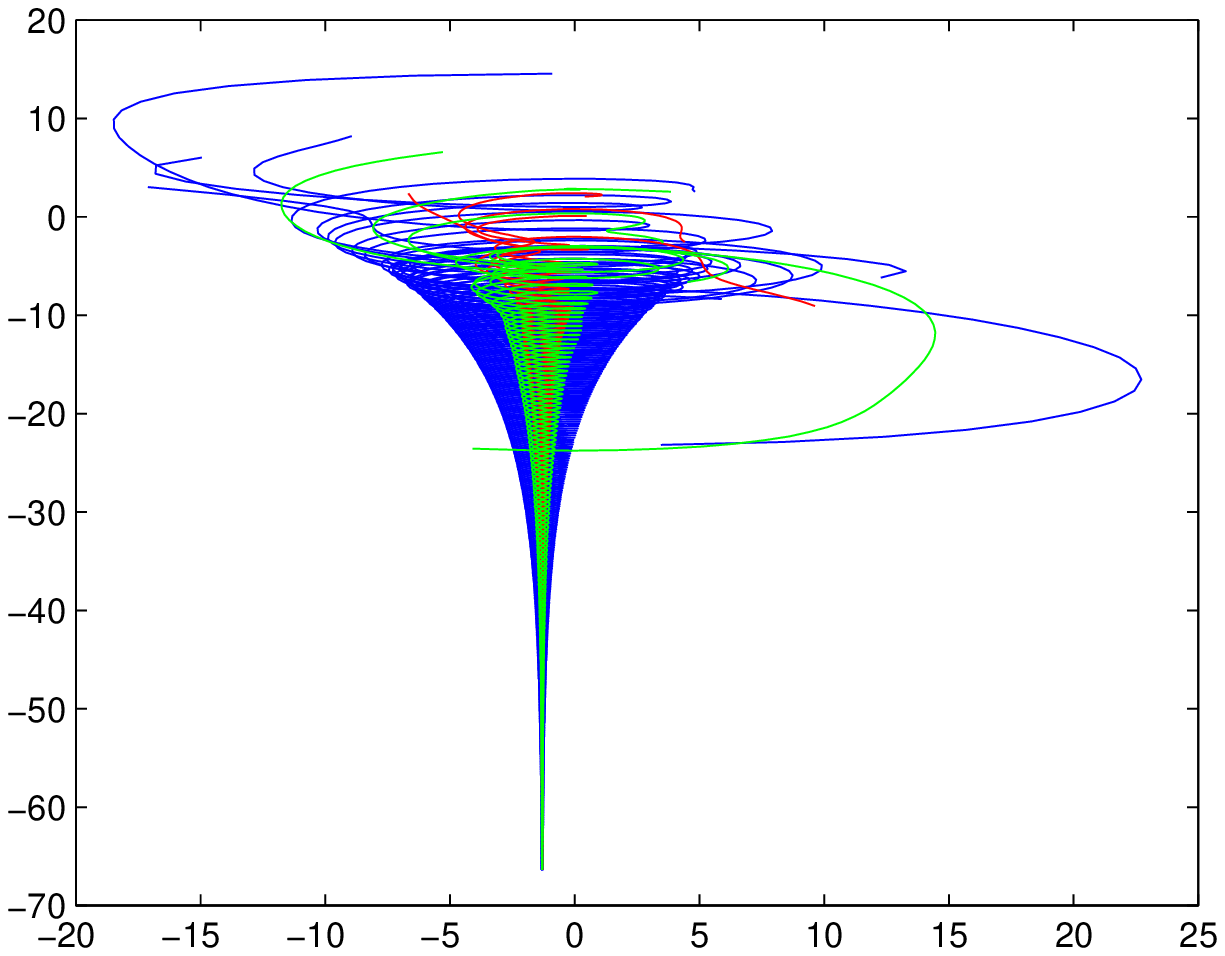}
\caption{Trajectories of 20 agents on same phase plane.
Trajectories of $S_1,\ldots,S_{10}$ colored blue,
$S_{11},\ldots,S_{15}$ colored red, 
and $S_{16},\ldots,S_{20}$ colored green.}
\label{dodec_phase}
\endminipage\hfill
\minipage{0.4\textwidth}
\endminipage\hfill
\minipage{0.1\textwidth}
\endminipage\hfill
\end{figure}

Let us consider digraph shown on Fig.\ref{dodec} consisting of 20 agents $S_1,\ldots,S_{20}$ described in section~\ref{agents_descr}.
This  dodecahedron-like digraph have Leading Set consisting of dynamic nodes $S_1,\ldots,S_{10}$ which are connected in directed circle.
Let us choose $\nu=k_1=k_2=\ldots=k_{10}$ and $\mu=k_{11}=k_{12}=\ldots=k_{20}.$  According to Theorem~\ref{th_circle} 
gain $\nu$ should be chosen $\nu>0.5\cdot\cot^2(\pi/10)\approx 4.74.$ For faster convergence $\nu$ let us choose $\nu=5.5.$
Simulations show that $\mu$ can be chosen considerably less than $\nu.$ Let us choose $\mu=1,$ and let agents have different
initial conditions. Results of numerical simulation show that such nonidentical gain choice provides achievement
of consensus. All trajectories of 20 agents on same phase plane are shown on Fig.~\ref{dodec_phase}.

Numerical simulations also show that 
by choosing $\mu=\nu$ and applying Theorem~\ref{th_compl_spectr} for resulting Laplace matrix one can obtain 
lower bound approximation  $\mu=\nu>4.74.$

\section{Reference remarks}
Assumptions of Th.~\ref{th1} are relaxed in comparison with Th. 2 from~\cite{fradkov_junussov_cdc11}.
Proof of these results use coordinate tranformation as in~\cite{YoshiokaNamerikawa08}. Lemma~\ref{lemma0}
partially succeeds Th. 3 from~\cite{FaxMurray04}. Th.~\ref{th_compl_spectr} is a trigonometric 
form of Th. 1 from~\cite{yu_cao_2010} with different proof.
\section{Conclusions}
By means of Passification Method sufficient conditions of consensus with identical
and nonidentical gains are derived. Synchronous behaviour (consensus vector) is described, it can
be affected by nonidentical gains (nonidentity in actuation) within Leading Set. Gain asymptote
in growing cycle digraphs which have lowest communication cost for reaching average consensus and consisting of
double integrators is studied. 

It is rediscovered that cycle digraph connection with nonidentical actuation of nodes 
causes nonidentical impact on synchronous behaviour. Reachability of synchronization corresponds to positive scalar -- common gain.
By constructing boundary of sufficient gain region in 3 node digraph it is found that Jordan block of Laplace matrix 
(which affects transient dynamics) appears in a direction of extremal point. Comparison of dynamics is a subject to a future study. 
Geometrical interpretations which might be useful in applications and theory were developed.
\begin{bibliography}{99}
\bibliographystyle{abbrv} 

\end{bibliography}

\end{document}